\documentclass[twoside,reqno,11pt]{amsart}
\pdfoutput=1
\usepackage[utf8]{inputenc}
\usepackage{xcolor}
\usepackage{prettyref}
\usepackage{mathtools}
\usepackage{amsmath}
\usepackage{amsthm}
\usepackage{amssymb}

\makeatletter
\numberwithin{equation}{section}
\numberwithin{figure}{section}

\theoremstyle{plain}
\newtheorem{thm}{\protect\theoremname}
\theoremstyle{plain}

\theoremstyle{plain}
\newtheorem{prop}[thm]{\protect\propositionname}
\theoremstyle{definition}

\theoremstyle{plain}
\newtheorem{lem}[thm]{\protect\lemmaname}
\theoremstyle{remark}
\newtheorem{remark}{Remark}

\@ifundefined{date}{}{\date{}}

\usepackage{textgreek}
\usepackage{graphicx}
\usepackage{amsfonts}
\usepackage{cases}
\usepackage{tikz-cd}
\usepackage{cancel}

\usepackage[framemethod=TikZ]{mdframed}

\def\XXint#1#2#3{{\setbox0=\hbox{$#1{#2#3}{\int}$ }
\vcenter{\hbox{$#2#3$ }}\kern-.6\wd0}}

\usepackage{accents}

\usepackage[unicode=true,pdfusetitle,
 bookmarks=true,bookmarksnumbered=true,bookmarksopen=false,
 breaklinks=false,pdfborder={0 0 1},backref=false]
 {hyperref}
\usepackage{cleveref}

\providetoggle{nomsort}
\settoggle{nomsort}{true} 

\makeatletter
\iftoggle{nomsort}{%
    \let\old@@@nomenclature=\@@@nomenclature        
        \newcounter{@nomcount} \setcounter{@nomcount}{0}%
        \renewcommand\the@nomcount{\two@digits{\value{@nomcount}}}
        \def\@@@nomenclature[#1]#2#3{
          \addtocounter{@nomcount}{1}%
        \def\@tempa{#2}\def\@tempb{#3}%
          \protected@write\@nomenclaturefile{}%
          {\string\nomenclatureentry{\the@nomcount\nom@verb\@tempa @[{\nom@verb\@tempa}]%
          \begingroup\nom@verb\@tempb\protect\nomeqref{\theequation}%
          |nompageref}{\thepage}}%
          \endgroup
          \@esphack}%
      }{}
\makeatother

\usepackage{tensind}
\tensordelimiter{?}

\hypersetup{citecolor=purple,linkcolor=blue}

\newcommand\restr[2]{{
  \left.\kern-\nulldelimiterspace 
  #1 
  \vphantom{\big|} 
  \right|_{#2} 
  }}

\usepackage{extarrows}

\usepackage{enumerate}\usepackage{color}\usepackage[margin=1.5in]{geometry}

\let\div\relax
\DeclareMathOperator\div{div}

\DeclareMathOperator\curl{curl}
\DeclareMathOperator\Id{Id}
\DeclareMathOperator\supp{supp}

\newcommand{\dd}[0]{\partial}
\newcommand{\grad}[0]{\nabla}

\newcommand{\T}[0]{\mathbb T^d}

\theoremstyle{plain}

\newcommand{\eq}[1]{\begin{align}{#1}\end{align}}
\newcommand{\eqn}[1]{\begin{align*}{#1}\end{align*}}
\newcommand{\Rd}[0]{{\mathbb R^2}}
\newcommand{\Td}[0]{{\mathbb T^2}}
\newcommand{\Ig}[2]{I^{g;#1}_{#2}}

\makeatother

\providecommand{\conjecturename}{Conjecture}
\providecommand{\definitionname}{Definition}
\providecommand{\lemmaname}{Lemma}
\providecommand{\propositionname}{Proposition}
\providecommand{\theoremname}{Theorem}

\begin{document}
\global\long\def\divop{\operatorname{div}}%
\global\long\def\supp{\operatorname{supp}}%
\global\long\def\curl{\operatorname{curl}}%

\global\long\def\grad{\nabla}%
\global\long\def\otherv{{u}}

\global\long\def\dd{\partial}%

\global\long\def\Ig{I_{\#2}^{g;\#1}}%

\global\long\def\T{\mathbb{T}^{2}}%

\title[Non-uniqueness for forced SQG]{Non-uniqueness up to the Onsager threshold for the forced SQG equation}

\author{Aynur Bulut}
\address{Louisiana State University, 303 Lockett Hall, Baton Rouge, LA 70803}
\email{aynurbulut@lsu.edu}
\author{Manh Khang Huynh}
\address{Georgia Institute of Technology Department of Mathematics, Atlanta, GA 30332}
\email{mhuynh41@gatech.edu}
\author{Stan Palasek}
\address{School of Mathematics, Institute for Advanced Study, Princeton, NJ 08540}
\email{palasek@ias.edu}

\begin{abstract}
We establish new non-uniqueness results for the forced inviscid surface quasi-geostrophic equation, via an alternating formulation of convex integration techniques.  Our results imply non-uniquenesss in the class of weak solutions with $|\grad|^{-1}\theta\in C_tC_x^\alpha$, for any $\alpha<1$.
\end{abstract}

\maketitle

\parskip0.1in

\allowdisplaybreaks

\section{Introduction}

Let $\Td$ denote the 2-torus $\Td\coloneqq (\mathbb R/2\pi\mathbb Z)^2$, which we identify with $[-\pi,\pi)^2$, and consider the two-dimensional surface quasi-geostrophic (SQG) equation
\begin{align} \label{sqg}
\left\lbrace\begin{array}{rl}\dd_t\theta+u\cdot\grad\theta&=f\\
u&=\Lambda^{-1}\grad^\perp\theta
\end{array}\right.
\end{align}
where $\theta:[0,\infty) \times\T\to\mathbb R$ is an unknown scalar field, $f:[0,\infty)\times\Td\to\mathbb R$ is an external force, $\grad^\perp\coloneqq(-\dd_2,\dd_1)^t$, and $\Lambda\coloneqq\sqrt{-\Delta}$.  The SQG equation is an important model in the physics of geophysical fluids; for a survey of the physical relevance of this system, see for instance \cite{salmon1998lectures}. 

From a mathematical standpoint, as first observed in \cite{constantin1994formation}, the equation \eqref{sqg} has significant interest as a model of potential singularity formulation for the Euler equations.  Indeed, the vector field $\omega=\grad^\perp\theta$ solves a two-dimensional analogue of the vorticity formulation of the three-dimensional Euler system.  Moreover, as in the case of the full three-dimensional Euler equations, global regularity of smooth finite energy solutions of \eqref{sqg} on $\T$ with smooth forcing remains an important open problem.  On the other hand, local well-posedness results (i.e. local in time existence and uniqueness of solutions) are known in the spaces $C^{1,\alpha}$ ($0<\alpha<1$) and $H^s$ ($s>2$); we refer the reader to \cite{constantin1994formation, majda2012compressible} for further discussion.  The question of whether these spaces are sharp for the local theory partially motivates our present work.

In particular, we focus on the question of whether weak solutions of the system \eqref{sqg} with the same data and force are uniquely determined.  To prepare for the statement and discussion of our main result, we begin with a brief review of some results from the literature in the unforced and forced settings.

\begin{enumerate}
\item[(i)] For the unforced system, corresponding to $f\equiv 0$, the question of uniqueness of weak solutions was answered in the negative for a certain class of low-regularity solutions by Buckmaster, Shkoller, and Vicol in \cite{buckmasterNonuniquenessWeakSolutions2019b}, based on the convex integration framework for fluid equations introduced by De Lellis and Sz\'ekelyhidi (see, e.g. \cite{de2013dissipative}) applied to a ``momentum'' formulation of the equation (see \eqref{eq:forced_sqg} below).  

Indeed, in \cite{buckmasterNonuniquenessWeakSolutions2019b} non-uniqueness is obtained by showing a stronger result: the construction of solutions with prescribed profile of the Hamiltonian $\left\Vert \theta(t)\right\Vert_{\dot{H}_{x}^{-1/2}}$ taken as any smooth compactly supported function.  This is closely related not only to non-uniqueness, but also to the problem of identifying the {\it Onsager threshold}, the regularity at which solutions to SQG start to conserve energy.  A similar Onsager-type theorem up to the same regularity was later obtained by Isett and Ma in \cite{isettma}, using different methods to cancel the stress in the original active scalar formulation \eqref{sqg}; see also \cite{mathesis}.  

The constructions in \cite{buckmasterNonuniquenessWeakSolutions2019b}, \cite{isettma} and \cite{mathesis} produce examples of non-unique solutions with regularity $\Lambda^{-1}\theta\in C_{t}^{\sigma}C_{x}^{\beta}$ where $\beta\in\left(\frac{1}{2},\frac{4}{5}\right)$ and $\sigma<\frac{\beta}{2-\beta}$, which falls short of the thresholds above which energy conservation and uniqueness are known to hold; these thresholds are $\beta=1$ and $\beta=2$, respectively.

\item[(ii)] In the forced setting, it is natural to ask whether one can construct forcing terms (of sufficient regularity) for which non-uniqueness holds beyond the regularity treated in \cite{buckmasterNonuniquenessWeakSolutions2019b}, \cite{isettma}, and \cite{mathesis}.  In our main result, stated below, we show that this is indeed the case.  

Our approach is motivated by \cite{bulutConvexIntegrationOnsager2023}, where the authors of the present paper developed an {\it alternating convex integration scheme} to show non-uniqueness for solutions to the forced Euler equations with regularity beyond the regularity for which energy is conserved (that is, beyond the Onsager threshold).  One goal of the present work is to show that the same strategy can be adapted to other fluid equations with forcing terms, and in particular in the case of the forced SQG equations.  In this context, we remark that there is also work on the existence of nontrivial steady solutions of SQG initiated by Cheng, Kwon, and Li \cite{cheng2021non} and improved in the forced case by Dai and Peng \cite{dai2023non} using a similar strategy as that in \cite{bulutConvexIntegrationOnsager2023}.
\end{enumerate}

We now state our main theorem, in which we establish non-uniqueness of weak solutions in a class where $\theta$ has barely negative regularity.  As in \cite{buckmasterNonuniquenessWeakSolutions2019b}, we consider the initial value problem for the forced SQG equations in their momentum formulation, given by 
\begin{align}
\left\lbrace\begin{array}{rl}\partial_{t}v+\Lambda v\cdot\nabla v-\left(\nabla v\right)^{T}\cdot\Lambda v+\nabla p & =\divop F\label{eq:forced_sqg}\\
\divop v & =0,\end{array}\right.
\end{align}
with initial data $v(0,x)=v_0(x)$ on $\mathbb{R}\times\mathbb{T}^{2}$, which has the effect of making the problem more closely resemble the Euler equations in velocity form (indeed, this choice is somewhat arbitrary, and we believe the methods in the present work could also be implemented in the context of the scalar equation, as in \cite{isettma}).

\begin{thm}[Non-uniqueness for forced SQG]
\label{thm:main_thm}  For any $\alpha\in(0,1)$, there exist distinct $\otherv,v\in C_{t}C_{x}^{\alpha}([0,\infty)\times\Td\to\mathbb R^2)$ and $F\in C_{t}C_{x}^{\alpha}([0,\infty)\times\Td\to\mathbb R^{2\times2})$ such that $\otherv(0,x)=v(0,x)$ and $\left(v,F\right)$ and $\left(\otherv,F\right)$ both solve \eqref{eq:forced_sqg}.
\end{thm}

Our use of the momentum formulation deserves some comment.  As in \cite{buckmasterNonuniquenessWeakSolutions2019b}, this formulation leads to expressions which possess favorable nonlinear interactions for convex integration (indeed, this was one of the main points in \cite{buckmasterNonuniquenessWeakSolutions2019b}, and is used to overcome the difficulty of the odd Fourier multiplier present in \eqref{sqg}; note that one can easily see that if the SQG nonlinearity is perturbed by a high-frequency flow in the usual way, the high-high-low interaction vanishes, and thus cannot be used to cancel with the Reynolds stress (see \cite{isettHolderContinuousSolutions2015})).

Note that one can recover solutions to the scalar SQG \eqref{sqg} from solutions to the system \eqref{eq:forced_sqg} by simply applying $-\grad^\perp\cdot$ to both sides; thus, comparing the variables to \eqref{eq:forced_sqg}, we have the correspondences $\theta=-\grad^\perp\cdot v$ and $f=-\grad^\perp\cdot\div F$.  In particular, solutions of \eqref{eq:forced_sqg} correspond to zero-average solutions of \eqref{sqg}, and the H\"older regularity of $u$, $v$, and $F$ in the statement of Theorem \ref{thm:main_thm} can be translated into the conclusion that the corresponding fields in the scalar SQG \eqref{sqg} can be made to satisfy $\theta_1,\theta_2\in \dot B^{-\epsilon}_{\infty,\infty}$ and $f\in \dot B^{-1-\epsilon}_{\infty,\infty}$ for any $\epsilon>0$.  

The proof of Theorem $\ref{thm:main_thm}$ is in the spirit of the alternating convex integration method we introduced in \cite{bulutConvexIntegrationOnsager2023}, and is based on a relaxation of \eqref{eq:forced_sqg} to a coupled forced SQG-Reynolds system,
\begin{align}\label{relaxedsystem}
\left\lbrace\begin{array}{rl}
\partial_{t}\otherv_{q}+\Lambda\otherv_{q}\cdot\nabla\otherv_{q}-\left(\nabla\otherv_{q}\right)^{T}\cdot\Lambda\otherv_{q}+\nabla\pi_{q} & =\divop F_{q}\\
\partial_{t}v_{q}+\Lambda v_{q}\cdot\nabla v_{q}-\left(\nabla v_{q}\right)^{T}\cdot\Lambda v_{q}+\nabla p_{q} & =\divop F_{q}+\divop R_{q}\\
\div u_q=\div v_q&=0,\end{array}\right.
\end{align}
for which we will construct solutions.  We refer to the first system as the {\it inactive} system, while the second system will be called the {\it active} system.  Roughly speaking, the construction proceeds by iteratively perturbing the active system to cancel the Reynolds stress, and then absorbing portions of the resulting error terms into the force (these steps are performed in alternate stages of the iteration).

When comparing our alternating SQG scheme to other instances where convex integration techniques are applied, a further obstacle arises from the lack of a ``gluing step'' which is a standard part of the Euler literature, following Isett \cite{isettProofOnsagerConjecture2018}. One notices that since the structure of the SQG nonlinearity does not match the form of the error $\div R_q$ (as it does in the case of Euler), it is not clear how to incorporate the error produced by the gluing into the Reynolds stress.\footnote{Expressing the SQG nonlinearity as the divergence of a Reynolds stress is possible using microlocal analysis (as shown in \cite{buckmasterNonuniquenessWeakSolutions2019b}) when the perturbation lives purely at high frequencies; unfortunately this is not the case with the perturbation caused by gluing.} Without gluing, there is no a priori way to obtain estimates on the material derivatives of the Reynolds stress (which are needed to estimate the transport error). The obvious remedy is to incorporate material derivative estimates into the iterative proposition; unfortunately this fails as well because after correcting, say, $v_q$, one has good estimates on the material derivative with respect to the $v_{q}$ flow, but not the $u_q$ one which would be needed for the next iteration step.

In the present work, we resolve this issue by keeping track of estimates on the Eulerian time derivative $\dd_tR_q$ in the iteration itself. Whereas in the Euler case this would yield far from optimal bounds, for SQG one can still obtain sharp estimates. Essentially, as long as the convecting flow (in this case $\Lambda v_q$; in Euler's case $v_q)$ measured in $L^\infty$ is subcritical during the iteration, the terms in the material derivative can be estimated separately. We remark that this use of $\dd_t$ in place of the material derivative simplifies the argument in some places, avoiding the need for some terms involving commutators of material derivatives and Fourier multipliers. 

We conclude this introduction with a few further remarks on our results.

\begin{remark}
The regularity of the constructed non-unique solutions in the statement of Theorem \ref{thm:main_thm} is sharp with respect to the Onsager threshold in the sense that as $\alpha\uparrow1$, the spaces containing the velocities and force potential approach $C_tC_x^{1+}$, and in this space energy conservation is known. On the other hand, there is still a gap between this result and the best known threshold for uniqueness, which is $u,v\in C_x^{2+}$.  Moreover, by well-known heuristic arguments, one should expect that the condition $\alpha<1$ in the statement of Theorem~\ref{thm:main_thm} is the best that can be obtained with current convex integration methods.  Indeed, one can compare this to the analogous threshold $\alpha<\frac12$ for the Euler equations, obtained in \cite{buckmaster2023intermittent} and \cite{bulutConvexIntegrationOnsager2023}; roughly speaking, one can argue that in order for the scheme to converge in a space with $\alpha$ derivatives, one must cancel out a stress $R$ living at frequencies $\lesssim\lambda$ with amplitude $A$ using a perturbation with amplitude $A^\frac12$ and frequency $\lambda'\gg\lambda$. If the nonlinearity of the system is, say, $n$-linear with $m$ derivatives, then one should have $A\sim(\lambda')^{-n\alpha}$. In order for the iteration to close, one needs
\eqn{
\lambda^m(\lambda')^{-m-n\alpha}\lesssim(\lambda'')^{-n\alpha}.
}
Rearranging, one finds that in order to have $\lambda$ increase rapidly, we are led to the condition $\alpha\leq m/n$. Thus in the case of the momentum\footnote{The scalar SQG, treated in \cite{isettma}, does not obey this heuristic because the cancellation of the stress by the nonlinearity is not as direct.} formulation of SQG, the expected threshold is $1$.
\end{remark}

\begin{remark}
It seems likely that ideas related to our proof of Theorem~\ref{thm:main_thm} can be used to prove an Onsager theorem in the sense that, up to the threshold $\theta\in C^{0-}$ (or $v\in C^{1-}$), there exist weak solutions of SQG that fail to conserve the ``Hamiltonian minus work'' quantity
\eqn{
\frac12\int_{\Td}|\Lambda^{-\frac12}\theta(t)|^2-\int_0^t\int_\Td\Lambda^{-1}\theta f
}
for some force $f$. Moreover this regularity class would be sharp; see, e.g.  \cite{isettHolderContinuousSolutions2015}. This is an interesting question that we do not pursue here.
\end{remark}

\begin{remark}
It is remarked in \cite{buckmasterNonuniquenessWeakSolutions2019b} that the issues obstructing convex integration from reaching the maximum regularity are similar for SQG and the 2D Euler equations---namely, that Beltrami flows are not particularly well-suited, in particular due to their poor interaction with the Lagrangian flow map. It is straightforward to see that the methods we develop in this paper (which reach the conjectured Onsager regularity $C^{0-}$ for SQG) can reach the threshold $C^{\frac13-}$ for the forced 2D Euler system. Recently this threshold has been reached for the unforced Euler system using different methods in \cite{giri2DOnsagerConjecture2023a}.
\end{remark}

\begin{remark}
In the Euler setting treated in \cite{bulutConvexIntegrationOnsager2023}, the force produced by the construction can be made small, and more generally can approximate any sufficiently smooth function.  It would be interesting to see whether the same is true in the forced SQG setting discussed in this paper.  We remark that it is still an open question as to how we can show non-uniqueness via convex integration above the Onsager exponent for the unforced Euler equations.
\end{remark}

\subsection*{Outline of the paper}

We now give a brief outline of the remainder of the paper.  In Section $2$, we introduce some notation and preliminaries, and introduce the main iterative construction underlying the proof of Theorem \ref{thm:main_thm}.  We formulate an iterative proposition, Proposition \ref{prop:iterative_proposition}, which will form the basis of this construction, and then give the proof of Theorem \ref{thm:main_thm} based on this statement.  The proof of Proposition \ref{prop:iterative_proposition} is then given in Section $3$.  We conclude the paper with two appendices, in which we respectively prove a useful geometric lemma used in our arguments, and establish several technical relationships enjoyed by the choice of parameters used in the iterative construction.

\subsection*{Acknowledgments}

The third author acknowledges support from the NSF under Grant No.\ DMS-1926686.

\section{Notation and Preliminaries}

For $\alpha\geq0$, we make use of the standard H\"older spaces $C^\alpha(\Td)$.  For brevity, we will write the norm as $\|\cdot\|_\alpha$. We often slightly abuse notation by writing $\|f\|_0$ for an $f=f(t,x)$ to denote the supremum over space and time.

Let us recall the definitions of several standard Fourier multipliers: the non-local derivative $\Lambda f(x)\coloneqq\sum_{k\in\mathbb Z^2}|k|\hat f(k)e^{ik\cdot x}$; the Leray projection $\mathbb Pf(x)\coloneqq\sum_{k\in\mathbb Z^2}(\Id-k\otimes k/|k|^2)\hat f(k)e^{ik\cdot x}$ onto the space of divergence-free vector fields; and the Littlewood-Paley projections $P_{\leq N}$ for each frequency scale $N>0$ (see, for instance, \cite{bahouriFourierAnalysisNonlinear2011} for the definition and basic properties).  We denote the material derivative along the flow $v$ by $D_{t,v}\coloneqq\dd_t+v\cdot\grad$.

We make use of the following asymptotic notation: $A\lesssim B$ or $A=O(B)$ indicate that there exists $c>0$ such that $|A|\leq cB$; moreover $c$ may not depend on any parameters except the final H\"older regularity of the construction ($\alpha$ and $\beta$). $A\ll B$ is a stronger statement by which we mean $|A|\leq cB$ where $c>0$ is a constant we can make arbitrarily small by taking the parameter $b$ close to $1$ or increasing $a$ (both the variables $a$ and $b$ are introduced in Section~\ref{framework}). Smallness of $c$ will not depend on any parameters except $\alpha$ and $\beta$, unless otherwise noted by adding a subscript, e.g., $A\ll_bB$. We also write $A\sim B$ to mean that $A\lesssim B\lesssim A$. Finally, for $a\in\mathbb R$, we use the notation $a-$ to indicate that the claim holds for $a-\epsilon$ for all sufficiently small $\epsilon>0$; moreover any implicit constants in the claim may depend on $\epsilon$.

\subsection{Convex integration framework and the proof of Theorem~\ref{thm:main_thm}}\label{framework}

We now outline our approach to the proof of Theorem~\ref{thm:main_thm}, which is based on the alternating scheme introduced by the authors in \cite{bulutConvexIntegrationOnsager2023} as an extension of the convex integration iteration framework; note that, in accordance with the statement of Theorem~\ref{thm:main_thm} and as in \cite{buckmasterNonuniquenessWeakSolutions2019b}, we consider the SQG system reformulated into its momentum version.  

We begin by fixing a parameter $\beta\in(\frac34,1)$ satisfying the condition $\beta>\max(\alpha,\frac34)$, and aim to show $u,v\in C_tC_x^\beta$ and $F\in C_tC_x^{2\beta-1}$.  This regularity condition will then ensure that $u$ and $v$ are as in the statement of the theorem.  We will use an iterative procedure to construct solutions to the relaxed SQG-Reynolds system \eqref{relaxedsystem}.

Fixing for the moment $a\gg1$ and $b\in\left(1,2\right)$, we introduce sequences of iteration parameters
\begin{align*}
\lambda_{q} & \coloneqq85\left\lceil a^{b^{q}}\right\rceil \\
\delta_{q}: & =\lambda_{q}^{-2\beta}.
\end{align*}

To keep track of estimates on the Eulerian time derivative $\dd_tR_q$ as described in the remarks following the statement of Theorem $\ref{thm:main_thm}$, we make a delicate choice of two different time scales: a mollification scale $\tau_{m,q+1}$ to avoid a loss-of-derivatives problem in time, and a cutoff scale $\tau_{c,q+1}$ to define Lagrangian coordinates on the appropriate time intervals. To get close to $\alpha=1$, one must take
\eqn{
\tau_{m,q+1}\ll\tau_{c,q+1}.
}
The particular values will be introduced Section~\ref{section-proof}. For the constraints on the parameters, see Appendix~\ref{sec:Parameter-calculus}.

We now introduce the main iterative proposition which underlies our arguments.  

\begin{prop}[Iterative proposition]
 \label{prop:iterative_proposition} Fix $q\in\mathbb{N}_{0},$ $\beta\in\left(\frac34,1\right)$, and suppose that we have $0<b-1\ll_{\beta}1$ and $1\ll_{b,\beta}a$. Let $u_q,v_q,R_q,F_q$ be smooth fields obeying \eqref{relaxedsystem} with suitable choices of the pressures. Suppose that there exists a constant $M>0$ (independent of $q$) such that we have the inductive estimates
\begin{align}
\left\Vert \otherv_{q}\right\Vert _{1} + \left\Vert \Lambda \otherv_{q}\right\Vert _{0}  \leq M\lambda_{q}\delta_{q}^{\frac{1}{2}},\quad 
\left\Vert v_{q}\right\Vert _{1} + \left\Vert \Lambda v_{q}\right\Vert _{0}   \leq M\lambda_{q-1}\delta_{q-1}^{\frac{1}{2}},\label{eq:ind_est0}
\end{align}
\begin{align}
\lambda_{q}^{-1}\left\Vert F_{q}\right\Vert _{1}+\left\Vert R_{q}\right\Vert _{0}+\lambda_{q}^{-1}\left\Vert R_{q}\right\Vert _{1} & \leq\epsilon\lambda_{q+1}\delta_{q+1},\label{eq:ind_est1}
\end{align}
and
\begin{align}
\left\Vert \partial_{t}F_{q}\right\Vert _{0}+\left\Vert \partial_{t}R_{q}\right\Vert _{0} & \leq\tau_{m,q}^{-1}\lambda_{q+1}\delta_{q+1}\label{eq:ind_est2}
\end{align}
where $\epsilon>0$ is an absolute constant chosen small enough to obey Lemma~\ref{lem:geometric_lemma}. Suppose also that
$v_{q}=P_{\lesssim\lambda_{q-1}}v_{q}$, $\otherv_{q}=P_{\lesssim\lambda_{q}}\otherv_{q}$, $F_{q}=P_{\lesssim\lambda_{q}}F_{q}$,
and $R_{q}=P_{\lesssim\lambda_{q}}R_{q}$, and that $v_{q},\otherv_{q},R_{q},F_{q}$
are identically zero on the time interval $\left[0,1-\sum_{i\leq q}\tau_{m,i}\right]$.

Then there exist $v_{q+1},\otherv_{q+1},R_{q+1}$, and $F_{q+1}$ such that
\begin{align}
\nonumber & F_{q+1}  =F_{q},\quad \otherv_{q+1}=\otherv_{q}, \\
 & \left\Vert v_{q+1}\right\Vert _{1} + \left\Vert \Lambda v_{q+1}\right\Vert _{0}  \leq M\lambda_{q+1}\delta_{q+1}^{\frac{1}{2}}, \quad \left\Vert v_{q+1}-v_{q}\right\Vert _{0}  \leq M\delta_{q+1}^{\frac{1}{2}}, \label{eq:v_q+1_C1}
\end{align} 
and
\begin{align}
\left\Vert R_{q+1}\right\Vert _{0}+\lambda_{q+1}^{-1}\left\Vert R_{q+1}\right\Vert _{1}  & \leq\epsilon\lambda_{q+2}\delta_{q+2}, \label{eq:R_q+1_C1} \\
\left\Vert \partial_{t}R_{q+1}\right\Vert _{0}  & \leq\tau_{m,q+1}^{-1}\lambda_{q+2}\delta_{q+2}, \label{eq:R_q+1_time}
\end{align}
and such that $v_{q+1}=P_{\lesssim\lambda_{q+1}}v_{q+1}$, $R_{q+1}=P_{\lesssim\lambda_{q+1}}R_{q+1}$.

Moreover, one can ensure that the spatial frequencies of $v_{q+1}-v_{q}$ are supported
near $$\left\{ k\lambda_{q+1}:k\in\Omega_{0}^{q\text{ mod }2}\cup\Omega_{1}^{q\text{ mod }2}\right\},$$
where $\Omega_{j}^{i}\subseteq S^{1}$ are mutually disjoint sets 
as defined in \prettyref{lem:geometric_lemma}, and that $v_{q+1}$ is
identically zero on the time interval $\left[0,1-\sum_{i\leq q+1}\tau_{m,i}\right]$.
\end{prop}

We defer the proof of Proposition~\ref{prop:iterative_proposition} to the next section.  Indeed, we now show that the statement of Proposition~\ref{prop:iterative_proposition} leads to a straightforward proof of \prettyref{thm:main_thm}.  

\begin{proof}[Proof of \prettyref{thm:main_thm}]
Suppose that Proposition~\ref{prop:iterative_proposition} holds, and let $V(t,x)$ be a given nonzero, smooth, time-dependent, divergence-free, mean-zero vector field such that $\supp_{t}V$ is a compact subset of $\left(1,\infty\right)$
and such that $V=P_{\leq1}V$. 

Set $v_{0}=V$ and $\otherv_{0}=-V$, and note that $F_{0}$ and $R_{0}$ can be determined from \eqref{relaxedsystem} by setting the pressure to zero and using the anti-divergence operator $\mathcal{B}$ defined in \eqref{eq:anti_divergence_operator}; this uses the fact that the SQG nonlinearity has zero mean, and thus the operator $\mathcal B$ acts as a right inverse of $\div$.

Fixing $\zeta>0$, and rescaling, for instance as in \cite{buckmasterOnsagerConjectureAdmissible2017},
\begin{align*}
v_{0}^{\zeta}\left(t,x\right) & {\,\coloneqq\,}\zeta v_{0}\left(\zeta t,x\right),\quad F_{0}^{\zeta}{\,\coloneqq\,}\zeta^{2}F_{0}\left(\zeta t,x\right)\\
\otherv_{0}^{\zeta}\left(t,x\right) & {\,\coloneqq\,}\zeta\otherv_{0}\left(\zeta t,x\right),\quad R_{0}^{\zeta}{\,\coloneqq\,}\zeta^{2}R_{0}\left(\zeta t,x\right)
\end{align*}
the inductive estimates required in \eqref{eq:ind_est0}--\eqref{eq:ind_est2} for $q=0$ are enforced upon choosing $\zeta>0$ sufficiently small. We remove the $\zeta$ from the notation and take $v_0,u_0,R_0,F_0$ to initialize the iteration, where the pressures chosen appropriately.

To construct subsequent iterates, we proceed as in \cite{bulutConvexIntegrationOnsager2023}. Suppose $q\geq0$ is even. Applying Proposition~\ref{prop:iterative_proposition} to $v_q,u_q,R_q,F_q$ yields a new solution tuple which we call $v_{q+1},u_{q+1},\tilde R_{q+1},\tilde F_{q+1}$. To effectively ``move'' the Reynolds stress from the $v$  system to the $u$ system, we redefine the Reynolds stress and pressure as
\eqn{
R_{q+1}\coloneqq -\tilde R_{q+1},\quad F_{q+1}\coloneqq \tilde F_{q+1}+\tilde R_{q+1}.
}
These new fields solve the modified version of \eqref{relaxedsystem},
\begin{equation*}\begin{aligned}
\partial_{t}\otherv_{q+1}+\Lambda\otherv_{q+1}\cdot\nabla\otherv_{q+1}-\left(\nabla\otherv_{q+1}\right)^{T}\cdot\Lambda\otherv_{q+1}+\nabla\pi_{q+1} & =\divop F_{q+1}+\divop R_{q+1}\\
\partial_{t}v_{q+1}+\Lambda v_{q+1}\cdot\nabla v_{q+1}-\left(\nabla v_{q+1}\right)^{T}\cdot\Lambda v_{q+1}+\nabla p_{q+1} & =\divop F_{q+1}\\
\div u_{q+1}=\div v_{q+1}&=0,
\end{aligned}\end{equation*}
One observes that $v_{q+1},u_{q+1},R_{q+1},F_{q+1}$ obey the hypotheses of Proposition~\ref{prop:iterative_proposition} \emph{with the roles of $v$ and $u$ reversed} and thus applies it to obtain $v_{q+2},u_{q+2},R_{q+2},F_{q+2}$. This scheme is iterated ad infinitum.

One computes that this process gives the estimates
\begin{align*}
\left\Vert v_{q+1}-v_{q}\right\Vert _{\beta-}+\left\Vert u_{q+1}-u_{q}\right\Vert _{\beta-} & \lesssim\lambda_{q+1}^{\beta-}\delta_{q+1}^{\frac{1}{2}}
\end{align*}
and
\begin{align*}
\|F_{q+1}-F_q\|_{2\beta-1-}\lesssim\left\Vert R_{q}\right\Vert _{\left(2\beta-1\right)-} & \lesssim\lambda_{q+1}^{2\beta-}\delta_{q+1},
\end{align*}
so that the sequences $\left(v_{q}\right)$, $\left(\otherv_{q}\right)$, and $\left(F_q\right)$ converge respectively to $v_{\infty} \in C_{t}^{0}C_{x}^{\beta-}$, $\otherv_{\infty}\in C_{t}^{0}C_{x}^{\beta-}$, and $F_\infty\in C_t^0C_x^{(2\beta-1)-}$. Similarly, $R_{q}$ converges to $0$ in $C_{t}^{0}C_{x}^{\left(2\beta-1\right)-}$. Clearly these notions of convergence suffice to conclude that $v_\infty$ and $u_\infty$ solve the momentum formulation of SQG with force $\div F_\infty$.

We observe that $v_{\infty}$ and $\otherv_{\infty}$ agree on the time interval $\left[0,1-\sum_{q\geq0}\tau_{m,q}\right]$, and thus these solutions have the same
initial data. Because $v_{0}\not\equiv\otherv_{0}$ and the spatial
frequencies of $v_{q+1}-v_{q}$ are supported near $\left\{ k\lambda_{q+1}:k\in\Omega_{0}^{0}\cup\Omega_{1}^{0}\right\} $,
while those of $\otherv_{q+1}-\otherv_{q}$ are supported near
$\left\{ k\lambda_{q+1}:k\in\Omega_{0}^{1}\cup\Omega_{1}^{1}\right\} $,
we can conclude that $v_{\infty}\not\equiv\otherv_{\infty}$.
\end{proof}

\section{Proof of Proposition~\ref{prop:iterative_proposition}}\label{section-proof}

We now turn our attention to the proof of Proposition~\ref{prop:iterative_proposition}.  Suppose that the inductive hypotheses given in the statement of the proposition are satisfied.  We begin with a standard mollification argument.

\subsection{Mollification}

We define a standard mollifier in time $\psi^{\tau_{m,q+1}}$ at scale
\[
\tau_{m,q+1}\coloneqq\left(\lambda_{q-1}\lambda_{q+1}\delta_{q-1}^{\frac{1}{2}}\right)^{-1},
\]
construct the mollified stress
\begin{align*}
R_{\ell} & =\psi^{\tau_{m,q+1}}*R_{q},
\end{align*}
and observe that the active system becomes
\begin{align*}
\partial_{t}v_{q}+\Lambda v_{q}\cdot\nabla v_{q}-\left(\nabla v_{q}\right)^{T}\cdot\Lambda v_{q}+\nabla p_{q} & =\divop F_{q}+\divop R_{\ell}+\divop(R_q-R_{\ell}).
\end{align*}
(Note that we do not mollify the velocity $v_q$ itself.) We view the term involving $R_q-R_{\ell}$ as a new error term, and observe that, by standard mollification estimates (see, for instance, \cite{contiHPrincipleRigidity12009}), we have the bound
\begin{align}
     \nonumber \left\Vert R_{q}-R_{\ell}\right\Vert _{0} &\lesssim\tau_{m,q+1}\left\Vert \partial_{t}R_{q}\right\Vert _{0}\lesssim \left(\lambda_{q-1}\lambda_{q+1}\delta_{q-1}^{\frac{1}{2}}\right)^{-1} \tau_{m,q}^{-1}\lambda_{q+1}\delta_{q+1}  \\   
  & =  \tau_{m,q}^{-1}{\lambda_{q-1}^{-1}}\delta_{q-1}^{-\frac{1}{2}}\delta_{q+1} \ll \epsilon \lambda_{q+2}\delta_{q+2}  \label{eq:rq_rl_C0}
\end{align}
where we have used \eqref{eq:2nd_inequality} and \eqref{eq:ind_est2}. 

We note the important detail that due to a gap in the exponents  (see Appendix~\ref{sec:Parameter-calculus}), the implicit constants in the estimates can be absorbed by increasing $a$.

We also have the material derivative estimate
\begin{align}
& \nonumber \left\Vert D_{t,\Lambda v_{q}}R_{q}\right\Vert _{0}+\left\Vert D_{t,\Lambda v_{q}}R_{\ell}\right\Vert _{0}  \lesssim\left\Vert \partial_{t}R_{q}\right\Vert _{0}+\left\Vert \Lambda v_{q}\right\Vert _{0}\left\Vert R_{q}\right\Vert _{1}\\
\nonumber & \lesssim\tau_{m,q}^{-1}\lambda_{q+1}\delta_{q+1}+\lambda_{q-1}\delta_{q-1}^{\frac{1}{2}}\cdot\lambda_{q}\lambda_{q+1}\delta_{q+1}\\
& =\lambda_{q+1}\delta_{q+1}\left(\tau_{m,q}^{-1}+\lambda_{q-1}\delta_{q-1}^{\frac{1}{2}}\lambda_{q}\right) \lesssim \lambda_{q+1}\delta_{q+1}\left(\tau_{m,q}^{-1}+\tau_{c,q+1}^{-1} \right)\label{eq:rq_rl_sep}
\end{align}
in which we used  \eqref{eq:1st_inequality}.
Let us also record the straightforward estimate
\begin{align}
\left\Vert \partial_{t}R_{\ell}\right\Vert _{0} & \lesssim\tau_{m,q+1}^{-1}\lambda_{q+1}\delta_{q+1}.\label{eq:time_deri_R}
\end{align}

\subsection{Perturbation}

We now introduce the main perturbation construction used in the convex integration step. For $k\in S^{1}$, we define the Beltrami plane waves
\begin{align*}
b_{k}(\xi)=ik^{\perp}e^{ik\cdot\xi}\in\mathbb{C}^2\quad\textrm{and}\quad c_{k}(\xi)=e^{ik\cdot\xi}\in\mathbb{C},
\end{align*}
so that $b_k(\xi)=\nabla^{\perp}c_{k}\left(\xi\right)$ and $c_k(\xi)=-\nabla^{\perp}\cdot b_{k}\left(\xi\right)$.  Note that $b_k$ and $c_k$ now satisfy
$$\Lambda b_{k}=i\left|k\right|k^{\perp}e^{ik\cdot\xi}=ik^{\perp}e^{ik\cdot\xi},$$
$\nabla\cdot b_{k}=0$, $\left|b_{k}\right|=1$, and $\left|c_{k}\right|=1$ for all $k\in S^1$.
We now fix a ``cutoff'' time scale 
\[
\tau_{c,q+1}\coloneqq\Big(\lambda_{q+1}\lambda_{q+2}\delta_{q+1}^{-\frac{1}{2}}\delta_{q+2}\Big)^{-1},
\]
and a sequence of times $(t_i)$, with $t_{i}\coloneqq i\tau_{c,q+1}$ for each $i=0,1,2,\ldots$. At each time $t_i$, we let $\Phi_{i}$ be the unique solution of the transport equation
\begin{align*}
D_{t,\Lambda v_{q}}\Phi_{i} & =0\\
\Phi_{i}(t_{i},x) & =x.
\end{align*}

Then, recalling \eqref{eq:1st_inequality}, we have the usual transport estimates
\begin{align}
\left\Vert \nabla\Phi_{i}(t)-\mathrm{Id}\right\Vert _{0} & \lesssim\left|t-t_{i}\right|\left\Vert \Lambda v_{q}\right\Vert _{1}\lesssim\tau_{c,q+1}\lambda_{q-1}^2\delta_{q-1}^{\frac{1}{2}}\ll1 \label{eq:transport_est_0}
\end{align}
and
\begin{align}
\left\Vert \nabla\Phi_{i}(t)\right\Vert _{N} & \lesssim\left|t-t_{i}\right|\left\Vert \Lambda v_{q}\right\Vert _{N+1}\lesssim\tau_{c,q+1}\lambda_{q-1}^{N+2}\delta_{q-1}^{\frac{1}{2}}\label{eq:transport_est}
\end{align}
for $N\geq 1$.  In what follows, we will also use the bound 
\begin{equation}
\left\Vert \partial_{t}\Phi_{i}\right\Vert _{0}=\left\Vert \Lambda v_{q}\cdot\nabla\Phi_{i}\right\Vert _{0}\lesssim\left\Vert \Lambda v_{q}\right\Vert _{0}\lesssim\lambda_{q-1}\delta_{q-1}^{\frac{1}{2}}.\label{eq:time_deri_phi}
\end{equation}


In order to properly cancel the Reynolds stress, we need the following geometric lemma which is a variant of Lemma~4.2 in \cite{buckmasterNonuniquenessWeakSolutions2019b}. The proof is deferred to Appendix~\ref{lemmaproof}.

\begin{lem}[Geometric lemma]
\label{lem:geometric_lemma} Let $B(\mathrm{Id},\epsilon)$ denote a ball around the identity in the space of symmetric $2\times 2$ matrices.  There exists $\epsilon>0$ such that there are \textbf{disjoint finite }subsets 
\[
\Omega_{0}^{0},\Omega_{1}^{0},\Omega_{0}^{1},\Omega_{1}^{1}\subseteq S^{1},
\]
and smooth positive functions 
\[
\gamma_{k}\in C^{\infty}\left(B_{\epsilon}\left(\mathrm{Id}\right)\right),\;\forall k\in\bigcup_{i,j}\Omega_{j}^{i}
\]
such that
\begin{enumerate}[(i)]
\item $85\Omega_{j}^{i}\subseteq\mathbb{Z}^{2}\;\forall i,j$
\item $\gamma_{k}=\gamma_{-k}$ and $-\Omega_{j}^{i}=\Omega_{j}^{i}\;\forall i,j$
\item We have 
\[
R=\frac{1}{2}\sum_{k\in\Omega_{j}^{i}}\gamma_{k}\left(R\right)^{2}\left(k^{\perp}\otimes k^{\perp}\right),\;\forall R\in B(\mathrm{Id},\epsilon),\forall i,j
\]
\item For $k,k'\in\Omega_{j}^{i}$, $k\neq-k'$, we have $\left|k+k'\right|>\frac12$.
\end{enumerate}
\end{lem}

Since $\lambda_q$ was defined as an integer multiple of $85$, property (i) implies that if $j\in\Omega_j^i$, then $\lambda_qk\in\mathbb Z^2$ for all $q\geq0$. This is essential so that the perturbation we construct is well-defined on $\Td$. 

Now, for each $i,j\in\mathbb{Z}$, define $\Omega_{j}^{i}\coloneqq\Omega^{i\mod 2}_{j\mod 2}$.  The motivating idea is that for each $j$, the set $\Omega_{j}^{i}$ corresponds to the truncation performed on a local time interval, and when gluing two adjacent intervals we should avoid interference. On the other hand, for a fixed $j$, the distinct projection sets $\Omega_{j}^{0}$ and $\Omega_j^1$ ensure that the sequences $\left(v_{q}\right)$ and $\left(\otherv_{q}\right)$ in the convex integration process have distinct spatial frequencies (implying the main non-uniqueness result in \prettyref{thm:main_thm}).

We next define the Fourier multiplier $\mathbb{P}_{q+1,k}=\mathbb{P}P_{\sim k\lambda_{q+1}}$
as a composition of Leray projection and frequency localization to
$\left|\frac{\xi}{\lambda_{q+1}}-k\right|=O(1)$, and define
a temporal partition of unity $(\chi_{i})\subset C_c^\infty(\mathbb{R})$, with $\chi_i\in C_{c}^{\infty}\left(\left(t_{i-1},t_{i+1}\right)\right)$ for each $i$, and such that 
\begin{align*}
\sum_{i}\chi_{i}^{2} & =1\text{ on }\supp_{t}\left(v_{q}\right)\cup\supp_{t}\left(F_{q}\right)\cup\supp_{t}\left(R_{q}\right)
\end{align*}
and
\begin{align*}
\sum_{i}\chi_{i}^{2} & =0\text{ near }0.
\end{align*}
Note that here and elsewhere, sums in the index $i$ are implicitly taken over the non-negative integers. With $(\gamma_{k})$ chosen as stated in Lemma~\ref{lem:geometric_lemma}, we then define 
\begin{align*}
a_{k}(t,x) & =\delta_{q+1}^{\frac{1}{2}}\gamma_{k}\left(\mathrm{Id}-\left(\lambda_{q+1}\delta_{q+1}\right)^{-1}R_{\ell}(t,x)\right)
\end{align*}
and
\begin{align*}
w_{q+1}\left(t,x\right) & =\sum_{i}\sum_{k\in\Omega_{i}^q}\chi_{i}\mathbb{P}_{q+1,k}\left(a_{k}\left(t,x\right)b_{k}\left(\lambda_{q+1}\Phi_{i}\left(t,x\right)\right)\right)\\
 & =\sum_{i}\sum_{k\in\Omega_{i}^q}\mathbb{P}_{q+1,k}\widetilde{w}_{q+1,i,k},
\end{align*}
where we have set $\widetilde{w}_{q+1,i,k}\coloneqq\chi_ia_kb_k(\lambda_{q+1}\Phi_i)$. Note, importantly, that 
\[\mathrm{Id}-\left(\lambda_{q+1}\delta_{q+1}\right)^{-1}R_{\ell}(t,x)\in B(\mathrm{Id},\epsilon)\]
by \eqref{eq:ind_est1} so Lemma~\ref{lem:geometric_lemma} applies; thus $\gamma_k$ is well-defined in the definition of $a_k$.

We make the observations that
\begin{align*}
D_{t,\Lambda v_{q}}w_{q+1}&=\sum_{i}\sum_{k\in\Omega_{i}^q}\mathbb{P}_{q+1,k}\left(\chi_{i}\left(D_{t,\Lambda v_{q}}a_{k}\right)b_{k}\left(\lambda_{q+1}\Phi_{i}\right)\right)\\
&\hspace{0.8in}+\mathbb{P}_{q+1,k}\left(\dd_t\chi_{i}a_{k}b_{k}\left(\lambda_{q+1}\Phi_{i}\right)\right)\\
&\hspace{0.8in}+\left[D_{t,\Lambda v_{q}},\mathbb{P}_{q+1,k}\right]\left(\chi_{i}a_{k}b_{k}\left(\lambda_{q+1}\Phi_{i}\right)\right)
\end{align*}
and, importantly,
\begin{align}
&\frac{R_{\ell}}{\lambda_{q+1}}+\frac12\sum_{k\in\Omega_{i}^q}a_{k}^{2}\left(k^{\perp}\otimes k^{\perp}\right)\nonumber\\
&\hspace{0.4in}=\frac{R_{\ell}}{\lambda_{q+1}}+\frac12\sum_{k\in\Omega_{i}^q}\delta_{q+1}\gamma_{k}^{2}\left(\mathrm{Id}-\left(\lambda_{q+1}\delta_{q+1}\right)^{-1}R_{\ell}\right)\left(k^{\perp}\otimes k^{\perp}\right)\nonumber\\
&\hspace{0.4in}=\frac{R_{\ell}}{\lambda_{q+1}}+\delta_{q+1}\left(\mathrm{Id}-\left(\lambda_{q+1}\delta_{q+1}\right)^{-1}R_{\ell}\right)\nonumber\\
&\hspace{0.4in}=\delta_{q+1}\Id.\label{Rcancellation}
\end{align}

Now, in view of \eqref{eq:ind_est1}, we have that, for each $N\geq 0$, 
\begin{equation}
\left\Vert a_{k}\right\Vert _{N}\lesssim\lambda_{q}^{N}\delta_{q+1}^{\frac{1}{2}}
\label{eq:a_est}
\end{equation}
and
\begin{align*}
\left\Vert \mathbb{P}_{q+1,k}\left(a_{k}b_{k}\left(\lambda_{q+1}\Phi_{i}\right)\right)\right\Vert _{N} & \lesssim \lambda_{q+1}^{N} \delta_{q+1}^{\frac{1}{2}}\left\Vert \gamma_{k}\left(\mathrm{Id}-\left(\lambda_{q+1}\delta_{q+1}\right)^{-1}R_{\ell}(t,x)\right)
e^{i k\cdot\lambda_{q+1}\Phi_{i}}\right\Vert _{0}\\
 & \lesssim\lambda_{q+1}^{N}\delta_{q+1}^{\frac{1}{2}}.
\end{align*}
where the implicit constants do not depend on $M$ from \eqref{eq:ind_est0}. 

We now estimate the norms of $w_{q+1}$, beginning with the $N=0$ case.  Note that 
\begin{align}
\left\Vert \widetilde{w}_{q+1,i,k}\right\Vert _{0} & \lesssim\delta_{q+1}^{\frac{1}{2}} \label{eq:w_tilde_est}
\end{align}
and, for $N\geq 0$, because $w_{q+1}= P_{\lesssim \lambda_{q+1}} w_{q+1}$   we have
\begin{align}
\left\Vert w_{q+1}\right\Vert _{N} + \left\Vert \Lambda^N w_{q+1}\right\Vert _{0}   & \lesssim\lambda_{q+1}^{N}\delta_{q+1}^{\frac{1}{2}},\nonumber
\end{align}
which proves  \eqref{eq:v_q+1_C1}. We remark that the presence of the constant $M$ in \eqref{eq:v_q+1_C1} comes from the implicit constants (independent of $q$) within our iteration process, and do not depend on the $M$ from the previous iterations in \eqref{eq:ind_est0}. This $M$ can not be absorbed by increasing $a$ as there was not a gap in the exponents to exploit (like in Appendix~\ref{sec:Parameter-calculus}).

By standard harmonic analysis \cite[Lemma A.8]{buckmasterNonuniquenessWeakSolutions2019b},
we have the commutator estimate
\[
\left\Vert \left[D_{t,\Lambda v_{q}},\mathbb{P}_{q+1,k}\right]\left(\widetilde{w}_{q+1,i,k}\right)\right\Vert _{0}\lesssim\left\Vert \nabla\Lambda v_{q}\right\Vert _{0}\left\Vert \widetilde{w}_{q+1,i,k}\right\Vert _{0}\lesssim\lambda_{q-1}^{2}\delta_{q-1}^{1/2}\delta_{q+1}^{\frac{1}{2}}
\]
Consequently, because of \eqref{eq:rq_rl_sep},
\begin{align*}
\left\Vert D_{t,\Lambda v_{q}}w_{q+1}\right\Vert _{0} & 
\lesssim  \left\Vert D_{t,\Lambda v_{q}}a_k\right\Vert _{0}+\tau_{c,q+1}^{-1}\delta_{q+1}^{\frac{1}{2}}+\lambda_{q-1}^{2}\delta_{q-1}^{1/2}\delta_{q+1}^{\frac{1}{2}}\\
 & \lesssim\delta_{q+1}^{\frac{1}{2}}\left(\tau_{m,q}^{-1}+ \tau_{c,q+1}^{-1} + \lambda_{q-1}\lambda_{q}\delta_{q-1}^{\frac{1}{2}}\right).
\end{align*}
Applying the parameter inequality \eqref{eq:1st_inequality}, we thus conclude
\begin{align}
\left\Vert D_{t,\Lambda v_{q}}w_{q+1}\right\Vert _{N} & \lesssim\lambda_{q+1}^{N}\delta_{q+1}^{\frac{1}{2}}\left(\tau_{m,q}^{-1}+\tau_{c,q+1}^{-1}\right) \label{eq:imprecise_d_tv_w}
\end{align}
for all $N\geq 0$, because $D_{t,\Lambda v_{q}}w_{q+1}$ has frequency support near the shell localized to $\lambda_{q+1}$.

Next, we note that for all $N\geq 0$, by \eqref{eq:time_deri_R} we have
\begin{align*}
\left\Vert \partial_{t}a_{k}\right\Vert _{N} & \lesssim\tau_{m,q+1}^{-1}\lambda_{q}^{N}\delta_{q+1}^{\frac{1}{2}},
\end{align*}
Applying the Leray projection to \eqref{relaxedsystem}, and using \eqref{eq:ind_est0}, \eqref{eq:ind_est1} and \eqref{eq:M_q+1},  we obtain 
\begin{align}\nonumber \left\Vert \partial_{t}\Lambda v_{q}\right\Vert _{0} & \lesssim \lambda_{q-1} \left( \left\Vert \nabla v_{q}\right\Vert _{0}\left\Vert \Lambda v_{q}\right\Vert _{0}+\left\Vert \divop\left(R_{q}\right)\right\Vert _{0}+\left\Vert \divop\left(F_{q}\right)\right\Vert _{0} \right) \\
 \nonumber  & \lesssim\lambda_{q-1}^{3}\delta_{q-1}+\lambda_{q}^{2}\lambda_{q+1}\delta_{q+1}\\
 & \lesssim\lambda_{q}^{2}\lambda_{q+1}\delta_{q+1} \ll\tau_{m,q+1}^{-1}\lambda_{q-1}\delta_{q-1}^{1/2} \label{eq:v_q_time}
\end{align}
As $\lambda_{q-1}\delta_{q-1}^{1/2}$ is the bound we have for $\left\Vert \Lambda v_{q}\right\Vert_0$, it effectively means that when the time derivative hits $\Lambda v_{q}$, we pay an extra factor of $\tau_{m,q+1}^{-1}$. Therefore,
\begin{align}
\nonumber \left\Vert \partial_{t}D_{t,\Lambda v_{q}}w_{q+1}\right\Vert _{N} & \lesssim
\lambda_{q+1}^{N} \left\Vert \partial_{t}D_{t,\Lambda v_{q}}w_{q+1}\right\Vert _{0} \\
& \nonumber \lesssim \lambda_{q+1}^{N}\delta_{q+1}^{\frac{1}{2}}\left(\tau_{m,q}^{-1}+\tau_{c,q+1}^{-1}\right)\cdot\left(\tau_{c,q+1}^{-1}+\tau_{m,q+1}^{-1}+\lambda_{q+1}\lambda_{q-1}\delta_{q-1}^{\frac{1}{2}}\right) \\
& \lesssim \lambda_{q+1}^{N}\delta_{q+1}^{\frac{1}{2}}\left(\tau_{m,q}^{-1}+\tau_{c,q+1}^{-1}\right)\cdot \tau_{m,q+1}^{-1}  \label{eq:w_q+1_transport_time} 
\end{align}
where we have used the bounds $\left|\partial_{t}\chi_{i}\right|\lesssim\tau_{c,q+1}^{-1}$, \eqref{eq:time_deri_R},  \eqref{eq:time_deri_phi} and \eqref{eq:v_q_time} (as the time derivative hits $\chi_i$, $a_k$,   $b_k(\lambda_{q+1}\Phi_i)$ and $\Lambda v_{q}$ respectively), and then apply the parameter inequality \eqref{eq:M_q+1}. In other words, $\tau_{m,q+1}^{-1}$ is the factor we pay for the time derivative.

\subsection{The new vector field and Reynolds stress}

At last, we may define the perturbed vector field
\[
v_{q+1}=v_{q}+w_{q+1}.
\]
In order to construct a suitable Reynolds stress so that $v_{q+1}$ solves the relaxed SQG with favorable estimates, one must define the anti-$\div$ operator $\mathcal{B}_0:C^\infty(\Td\to\Rd)\to C^\infty(\Td\to\mathbb R^{2\times2})$ by
\begin{align*}
\left(\mathcal{B}_0f\right)^{ij} & \coloneqq-\left(-\Delta\right)^{-1}\left(\partial_{i}f^{j}+\partial_{j}f^{i}\right),
\end{align*}
which takes values in the space of symmetric $2\times2$ matrices and satisfies
\begin{align*}
\left(\divop\mathcal{B}_0f\right)^{i} & =\partial_{j}\left(\mathcal{B}_0f\right)^{ij}=-\left(-\Delta\right)^{-1}\left(\partial_{i}\divop f\right)+f^{i}-|\Td|^{-1}\int_\Td f^i dx.
\end{align*}
This property is unchanged when the operator is pre-composed with the Leray projection,
\eq{
\mathcal B\coloneqq\mathcal B_0\mathbb P\label{eq:anti_divergence_operator}
}
Thus we have
\eq{
\div\mathcal Bf=\mathbb P f-|\Td|^{-1}\int_\Td f dx\label{antidiv}
}
for all smooth vector fields $f$. Moreover, $\mathcal B$ sends $\grad f$ to zero for any scalar field $f$.

Now we are in a position to define the new Reynolds stress $R_{q+1}$. After applying the identity
\eq{\label{identity}
\left(\nabla w_{q+1}\right)^{T}\cdot\Lambda v_{q}+\left(\nabla\Lambda v_{q}\right)^{T}\cdot w_{q+1}=\nabla\left(w_{q+1}\cdot\Lambda v_{q}\right),
}
one computes that $v_{q+1}$ obeys the system
\begin{align*}
&\partial_{t}v_{q+1}+\Lambda v_{q+1}\cdot\nabla v_{q+1}-\left(\nabla v_{q+1}\right)^{T}\cdot\Lambda v_{q+1}+\nabla p_{q+1}\\
&\hspace{0.8in}=\partial_{t}\left(v_{q}+w_{q+1}\right)+\Lambda\left(v_{q}+w_{q+1}\right)\cdot\nabla\left(v_{q}+w_{q+1}\right)\\
&\hspace{1.4in}-\left(\nabla\left(v_{q}+w_{q+1}\right)\right)^{T}\cdot\Lambda\left(v_{q}+w_{q+1}\right)+\nabla p_{q+1}\\
&\hspace{0.8in}=\divop F_{q}+\divop\left(R_{q}-R_{\ell}\right)+\divop R_{\mathrm{osc}}+\divop R_{\mathrm{tran}}\\
&\hspace{1.4in}+\divop R_{\mathrm{Nash}}-\nabla\left(w_{q+1}\cdot\Lambda v_{q}\right),
\end{align*}
where $R_{\mathrm{osc}}$, $R_{\mathrm{tran}}$, and $R_{\mathrm{Nash}}$ will be defined to satisfy
\begin{align*}
\div R_{\mathrm{osc}} & =  \divop R_{\ell}+\Lambda w_{q+1}\cdot\nabla w_{q+1}-\left(\nabla w_{q+1}\right)^{T}\cdot\Lambda w_{q+1} \mod\grad C^\infty \\
\div R_{\mathrm{tran}} & = \partial_{t}w_{q}+\Lambda v_{q}\cdot\nabla w_{q+1}\mod\grad C^\infty\\
\div R_{\mathrm{Nash}} & =  \Lambda w_{q+1}\cdot\nabla v_{q}-\left(\nabla v_{q}\right)^{T}\cdot\Lambda w_{q+1}+\left(\nabla\Lambda v_{q}\right)^{T}\cdot w_{q+1} \mod\grad C^\infty.
\end{align*}
There are places in the argument where extra irrotational vector fields appear leading to these identities being true only ``modulo $\grad C^\infty$'', such as the inversion of $\div$ by $\mathcal B$ (recall \eqref{antidiv}); and the low-frequency interactions in $R_{\mathrm{osc}}$ (see Section~\ref{oscerror}). The pressure is then modified to absorb them.

To proceed, we set
\begin{align*}
R_{q+1}&=\left(R_{q}-R_{\ell}\right)+R_{\mathrm{osc}}+R_{\mathrm{tran}}+R_{\mathrm{Nash}}.
\end{align*}

The rest of the proof is to show \eqref{eq:R_q+1_C1}  and \eqref{eq:R_q+1_time}.

The mollification error $R_{q}-R_{\ell}$ has already been estimated from \eqref{eq:ind_est1}, \eqref{eq:ind_est2} and \eqref{eq:rq_rl_C0}, and any implicit constants (which do not depend on $q$) can be absorbed by increasing $a$.

We will estimate the transport, Nash, and oscillation errors individually.

\subsubsection{Transport error}

We first estimate the transport error $R_{\mathrm{tran}}$ which is defined as
\eqn{
R_{\mathrm{tran}}\coloneqq\mathcal{B}D_{t,\Lambda v_{q}}w_{q+1}.
}
For this, we note that $w_{q+1}$ has frequency support near the shell
localized to $\lambda_{q+1}$ while $v_{q}$ has support near the shells up to $\lambda_{q-1}$,
so $\Lambda v_{q}\cdot\nabla w_{q+1}$ has support near $\lambda_{q+1}$.  Thanks to this, we can apply the stationary phase argument that has become standard in the convex integration literature: that applying a $-1$-order operator (in this case $\mathcal B$) to something frequency-localized in this way yields a gain of a factor $\lambda_{q+1}^{-1}$. By Bernstein's inequality, \eqref{eq:imprecise_d_tv_w}, and \eqref{eq:implication_un}, we therefore have
\begin{align*}
\left\Vert R_{\mathrm{tran}}\right\Vert _{0}&=\left\Vert \mathcal{B}D_{t,\Lambda v_{q}}w_{q+1}\right\Vert _{0}\\
&\lesssim\lambda_{q+1}^{-1}\left\Vert D_{t,\Lambda v_{q}}w_{q+1}\right\Vert _{0}\\
&\lesssim\lambda_{q+1}^{-1}\delta_{q+1}^{\frac{1}{2}}\left(\tau_{m,q}^{-1}+\tau_{c,q+1}^{-1}\right)\\
&\ll \epsilon \lambda_{q+2}\delta_{q+2}.
\end{align*}
We observe that $R_{\mathrm{tran}}= P_{\lesssim \lambda_{q+1}} R_{\mathrm{tran}}$ so the bound 
\[\left\Vert R_{\mathrm{tran}} \right \Vert_1 \ll \epsilon \lambda_{q+1} \lambda_{q+2}\delta_{q+2} \] 
immediately follows. Once again we note that the implicit constants can be absorbed by increasing $a$ due to gaps in the exponents.

Turning to the time derivative $\partial_tR_{\mathrm{tran}}$, note that there is a loss of derivatives in estimating $\Vert \partial_{t}R_{\mathrm{tran}}\Vert _{0}$ (due to the presence of two time derivatives on $w_{q+1}$).  Fortunately, as in \eqref{eq:w_q+1_transport_time}, another time derivative on $D_{t,\Lambda v_{q}}w_{q+1}$ incurs a factor of at most 
$\tau_{m,q+1}^{-1}$. Therefore, in view of \eqref{eq:M_q+1} and how we estimated $\left\Vert R_{\mathrm{tran}} \right \Vert_0$ above, we have 
\[
\left\Vert \partial_{t}R_{\mathrm{tran}}\right\Vert _{0}\ll  \tau_{m,q+1}^{-1}\lambda_{q+2}\delta_{q+2}.
\]
which proves \eqref{eq:R_q+1_time} for the transport error.
\subsubsection{Nash error}

We now estimate $R_{\mathrm{Nash}}$ which is defined as
\eqn{
R_{\mathrm{Nash}}=\mathcal{B}\left(\Lambda w_{q+1}\cdot\nabla v_{q}-\left(\nabla v_{q}\right)^{T}\cdot\Lambda w_{q+1}+\left(\nabla\Lambda v_{q}\right)^{T}\cdot w_{q+1}\right).
}
Here, the stationary phase argument noted above along with \eqref{eq:1st_inequality} gives the estimate
\begin{align*}
\left\Vert \mathcal{B}\left(\left(\nabla\Lambda v_{q}\right)^{T}\cdot w_{q+1}\right)\right\Vert _{0} & \lesssim\lambda_{q+1}^{-1}\left\Vert \left(\nabla\Lambda v_{q}\right)^{T}\cdot w_{q+1}\right\Vert _{0}\\
&\lesssim\lambda_{q+1}^{-1}\lambda_{q-1}^{2}\delta_{q-1}^{\frac{1}{2}}\delta_{q+1}^{\frac{1}{2}}\\
&\lesssim\lambda_{q+1}^{-1}\lambda_{q-1}\lambda_{q}\delta_{q-1}^{\frac{1}{2}}\delta_{q+1}^{\frac{1}{2}}\\
&\ll \epsilon \lambda_{q+2}\delta_{q+2}.
\end{align*}

For convenience, we introduce the notation
\begin{align*}
\psi_{q+1,j,k}\left(t,x\right) & =e^{i\lambda_{q+1}\left(\Phi_{j}\left(t,x\right)-x\right)\cdot k}\\
\psi_{q+1,j,k}\left(t,x\right)c_{k}\left(\lambda_{q+1}x\right) & =c_{k}\left(\lambda_{q+1}\Phi_{j}\left(t,x\right)\right)\\
\psi_{q+1,j,k}\left(t,x\right)b_{k}\left(\lambda_{q+1}x\right) & =b_{k}\left(\lambda_{q+1}\Phi_{j}\left(t,x\right)\right)
\end{align*}
(defining $\psi_{q+1,j,k}$ so that the operator $\nabla\Phi_{j}-\mathrm{Id}$ appears; this yields a crucial gain as seen in \eqref{eq:transport_est_0} and \eqref{eq:transport_est}), so one can re-express the perturbation as
\begin{align*}
w_{q+1}\left(t,x\right) & =\sum_{i}\sum_{k\in\Omega_{i}^q}\chi_{i}\mathbb{P}_{q+1,k}\left(a_{k}\psi_{q+1,i,k}b_{k}\left(\lambda_{q+1}x\right)\right).
\end{align*}
Therefore, recalling \eqref{eq:transport_est}, for $N\geq 1$ we have the bound
\begin{align*}
\left\Vert \psi_{q+1,j,k}\right\Vert _{N} & \lesssim_N \left(\lambda_{q+1}\left\Vert \nabla\Phi_{j}-\mathrm{Id}\right\Vert _{0}\right)^{N}+\lambda_{q+1}\left\Vert \nabla\Phi_{j}-\mathrm{Id}\right\Vert _{N-1}\nonumber \\
 & \lesssim\left(\lambda_{q+1}\tau_{c,q+1}\lambda_{q-1}^2\delta_{q-1}^{\frac{1}{2}}\right)^{N}+\lambda_{q+1}\tau_{c,q+1}\lambda_{q-1}^{N+1}\delta_{q-1}^{\frac{1}{2}}\\
 &\lesssim \lambda_q^N
\end{align*}
by \eqref{eq:1st_derived}. From here one concludes that for all $N\geq 1$,
\begin{align} \left\Vert a_{k}\psi_{q+1,j,k}\right\Vert _{N} & \lesssim_N \lambda_{q}^{N}\delta_{q+1}^{\frac{1}{2}},\label{eq:a_psi_est}
\end{align}
using \eqref{eq:a_est}.

We briefly remark that when the time derivative hits $ \psi_{q+1,j,k}$, we pay a factor of $\tau_{m,q+1}^{-1}$, just like when the time derivative hits $b_k(\lambda_{q+1}\Phi_j)$ in \eqref{eq:w_q+1_transport_time}.

To estimate the remaining two terms, note that
\begin{align*}
 & \left\Vert \mathcal{B}\left(\Lambda w_{q+1}\cdot\nabla v_{q}-\left(\nabla v_{q}\right)^{T}\cdot\Lambda w_{q+1}\right)\right\Vert _{0}=\left\Vert \mathcal{B}\left(\Lambda w_{q+1}^{\perp}\left(\nabla^{\perp}\cdot v_{q}\right)\right)\right\Vert _{0}.
\end{align*}
Moreover, observing that we have
\begin{align*}
b_{k}^{\perp}\left(\xi\right) & =\nabla^{\perp\perp}c_{k}\left(\xi\right)=-\nabla c_{k}\left(\xi\right),
\end{align*}
and thus
\begin{align*}
b_{k}^{\perp}(\lambda x) & =-\nabla c_{k}\left(\lambda x\right)=-\frac{1}{\lambda}\nabla\left(c_{k}\left(\lambda x\right)\right),
\end{align*}
it follows that we have the identity
\begin{align*}
\Lambda w_{q+1}^{\perp} & =\sum_{i}\sum_{k\in\Omega_{i}^q}\Lambda\mathbb{P}_{q+1,k}\left(\chi_{i}a_{k}\psi_{q+1,i,k}b_{k}^{\perp}\left(\lambda_{q+1}x\right)\right)\\
 & =-\frac{1}{\lambda_{q+1}}\sum_{i}\sum_{k\in\Omega_{i}^q}\Lambda\mathbb{P}_{q+1,k}\left(\chi_{i}a_{k}\psi_{q+1,i,k}\nabla\left(c_{k}\left(\lambda_{q+1}x\right)\right)\right)\\
 & =-\frac{1}{\lambda_{q+1}}\nabla\sum_{i}\sum_{k\in\Omega_{i}^q}\Lambda\mathbb{P}_{q+1,k}\left(\chi_{i}a_{k}\psi_{q+1,i,k}c_{k}\left(\lambda_{q+1}x\right)\right)\\
 & \phantom{=\frac{1}{\lambda_{q+1}}}+\frac{1}{\lambda_{q+1}}\sum_{i}\sum_{k\in\Omega_{i}^q}\Lambda\mathbb{P}_{q+1,k}\left(\nabla\left(\chi_{i}a_{k}\psi_{q+1,i,k}\right)c_{k}\left(\lambda_{q+1}x\right)\right).
\end{align*}
Now, recall that $\mathcal{B}\left(\nabla f\right)=0$ for any scalar $f$.  We therefore obtain (via an ``anti-div by parts'' argument),
\begin{align*}
&\mathcal{B}\left(\Lambda w_{q+1}^{\perp}\left(\nabla^{\perp}\cdot v_{q}\right)\right)\\
&\hspace{0.2in}=\frac{1}{\lambda_{q+1}}\sum_{i}\sum_{k\in\Omega_{i}^q}\mathcal{B}\left(\nabla\left(\nabla^{\perp}\cdot v_{q}\right)\cdot\Lambda\mathbb{P}_{q+1,k}\left(\chi_{i}a_{k}\psi_{q+1,i,k}c_{k}\left(\lambda_{q+1}x\right)\right)\right)\\
&\hspace{0.2in}\phantom{=\frac{1}{\lambda_{q+1}}}+\frac{1}{\lambda_{q+1}}\sum_{i}\sum_{k\in\Omega_{i}^q}\mathcal{B}\left(\left(\nabla^{\perp}\cdot v_{q}\right)\cdot\Lambda\mathbb{P}_{q+1,k}\left(\nabla\left(\chi_{i}a_{k}\psi_{q+1,i,k}\right)c_{k}\left(\lambda_{q+1}x\right)\right)\right).
\end{align*}
Using \eqref{eq:a_psi_est}, we therefore have 
\begin{align*}
&\left\Vert \mathcal{B}\left(\Lambda w_{q+1}^{\perp}\left(\nabla^{\perp}\cdot v_{q}\right)\right)\right\Vert _{0}\\
&\hspace{0.4in}\lesssim \frac{1}{\lambda_{q+1}^{2}}\sum_{i}\sum_{k\in\Omega_{i}^q}\left\Vert v_{q}\right\Vert _{2}\cdot\left\Vert \chi_{i}a_{k}\psi_{q+1,i,k}c_{k}\left(\lambda_{q+1}x\right)\right\Vert _{1}\\
&\hspace{0.4in}\phantom{=\frac{1}{\lambda_{q+1}}}+\frac{1}{\lambda_{q+1}^{2}}\sum_{i}\sum_{k\in\Omega_{i}^q}\left\Vert v_{q}\right\Vert _{1}\cdot\left\Vert \nabla\left(\chi_{i}a_{k}\psi_{q+1,i,k}\right)c_{k}\left(\lambda_{q+1}x\right)\right\Vert _{1}\\
&\hspace{0.4in}\lesssim \frac{1}{\lambda_{q+1}^{2}}\left(\lambda_{q-1}^{2}\delta_{q-1}^{\frac{1}{2}}\right)\cdot\delta_{q+1}^{\frac{1}{2}}\lambda_{q+1}+\frac{1}{\lambda_{q+1}^{2}}\left(\lambda_{q-1}\delta_{q-1}^{\frac{1}{2}}\right)\cdot\delta_{q+1}^{\frac{1}{2}}\left(\lambda_{q}\lambda_{q+1}\right)\\
&\hspace{0.4in}\lesssim  \frac{1}{\lambda_{q+1}}\lambda_{q-1}\lambda_{q}\delta_{q-1}^{\frac{1}{2}}\delta_{q+1}^{\frac{1}{2}}\\
&\hspace{0.4in}\ll \epsilon \lambda_{q+2}\delta_{q+2},
\end{align*}
where the last inequality follows from \eqref{eq:1st_inequality}. 

Collecting these estimates, and by arguing as above for the transport error, we obtain \eqref{eq:R_q+1_C1}
 and \eqref{eq:R_q+1_time} for the Nash error.
 
\subsubsection{Oscillation error}\label{oscerror}

We now turn to the construction and estimate of the oscillation error $R_{\mathrm{osc}}$. This is a very delicate part of the argument but fortunately the method presented in \cite{buckmasterNonuniquenessWeakSolutions2019b} is still effective in this setting. The general idea is to split into two terms 
\eqn{
R_{\mathrm{osc}}=R_{\mathrm{osc,high}}+R_{\mathrm{osc,low}}
}

We start with the observation that for $k,k'\in\Omega^i_{j}$ with $k\neq k'$, we have $\left|k+k'\right| \sim 1$. It follows that for general $k\in\Omega^q_{j},k'\in\Omega^{q}_{j'}$ with $k\neq-k'$ we have $\left|k+k'\right|\sim1$, and thus $\left(\Lambda\mathbb{P}_{q+1,k}\widetilde{w}_{q+1,i,k}\right)\cdot\left(\nabla\mathbb{P}_{q+1,k'}\widetilde{w}_{q+1,i,k'}\right)$ has frequency support on the shell localized to $\lambda_{q+1}$. These are the fortunate high-frequency parts (when $k\neq k'$). In view of this, we define
\eqn{
R_{\mathrm{osc,high}}\coloneqq\mathcal{B}P_{\sim\lambda_{q+1}}\sum_{i,i'}\sum_{k+k'\neq0}\Big(\Lambda\mathbb{P}_{q+1,k}\widetilde{w}_{q+1,i,k}\cdot\nabla\mathbb{P}_{q+1,k}\widetilde{w}_{q+1,i',k'}\\
-\left(\nabla\mathbb{P}_{q+1,k}\widetilde{w}_{q+1,i,k}\right)^{T}\cdot\Lambda\mathbb{P}_{q+1,k}\widetilde{w}_{q+1,i',k'}\Big).
}
We also define the low-frequency part which consists of the interactions where $k+k'=0$ (which implies that $k$ and $k'$ are in the same $\Omega_{i}^q$ for some $i$). Thus we have the decomposition
\eqn{
\divop R_{\mathrm{osc,low}} & = \divop R_{\ell}+\sum_{i,k}\mathfrak{I}_{i,k} \mod \grad C^{\infty}
}
where
\begin{align*}
\mathfrak{I}_{i,k} & =\frac{1}{2}\left(\Lambda\mathbb{P}_{q+1,k}\widetilde{w}_{q+1,i,k}\cdot\nabla\mathbb{P}_{q+1,k}\widetilde{w}_{q+1,i,-k}\right)\\
&\hspace{0.4in}+\frac{1}{2}\left(\Lambda\mathbb{P}_{q+1,k}\widetilde{w}_{q+1,i,-k}\cdot\nabla\mathbb{P}_{q+1,k}\widetilde{w}_{q+1,i,k}\right)\\
&\hspace{0.4in}-\frac{1}{2}\left(\nabla\mathbb{P}_{q+1,k}\widetilde{w}_{q+1,i,k}\right)^{T}\cdot\Lambda\mathbb{P}_{q+1,k}\widetilde{w}_{q+1,i,-k}\\
&\hspace{0.4in}-\frac{1}{2}\left(\nabla\mathbb{P}_{q+1,k}\widetilde{w}_{q+1,i,-k}\right)^{T}\cdot\Lambda\mathbb{P}_{q+1,k}\widetilde{w}_{q+1,i,k},
\end{align*}
which is a symmetric expression, i.e.  $\mathfrak{I}_{i,k}=\mathfrak{I}_{i,-k}$.

By the calculations in Section~5.4.3 in \cite{buckmasterNonuniquenessWeakSolutions2019b}, we have a further decomposition
\begin{align}
&\divop R_{\mathrm{osc,low}}=\divop\bigg(R_{\ell}+\sum_{i}\sum_{k\in\Omega_{i}^{q}} \Big( \frac{\lambda_{q+1}}{2}\chi_{i}^{2}\left(-k\otimes k\right)a_{k}^{2}\nonumber\\
&\hspace{2.6in}+\widetilde{\mathcal{Q}}_{i.k}^{\left(1\right)}+\mathcal{\widetilde{Q}}_{i.k}^{\left(2\right)} \Big) \bigg) \mod \grad C^{\infty}\label{osc_low_decomp}
\end{align}
where 
\begin{align*}
\left(\widetilde{\mathcal{Q}}_{i.k}^{\left(1\right)}\right)^{jl}\left(t,x\right) & =\frac{\chi_{i}^{2}}{2}\int_{0}^{1}dr\int_{0}^{1}d\overline{r}\int_{\mathbb{R}^{2}\times\mathbb{R}^{2}}dz_{1}dz_{2}\\
 & \;\left(\mathcal{K}_{k,r,\overline{r}}^{\left(1\right)}\right)^{jl}\left(x-z_{1},x-z_{2}\right)\cdot\nabla\left(a_{k}\psi_{q+1,i,k}\right)\left(z_{1}\right)a_{k}\psi_{q+1,i,-k}\left(z_{2}\right)\\
\left(\widetilde{\mathcal{Q}}_{i.k}^{\left(2\right)}\right)^{jl}\left(t,x\right) & =\frac{\chi_{i}^{2}}{2}\int_{0}^{1}dr\int_{0}^{1}d\overline{r}\int_{\mathbb{R}^{2}\times\mathbb{R}^{2}}dz_{1}dz_{2}\\
 & \;\left(\mathcal{K}_{k,r,\overline{r}}^{\left(2\right)}\right)^{jl}\left(x-z_{1},x-z_{2}\right)\cdot\nabla\left(a_{k}\psi_{q+1,i,-k}\right)\left(z_{2}\right)a_{k}\psi_{q+1,i,k}\left(z_{1}\right)
\end{align*}
in which $\mathcal{K}_{k,r,\overline{r}}^{(j)}$ are some kernels which satisfy the  estimates 
\begin{align}
\left\Vert \left(z_{1},z_{2}\right)^{a}\nabla_{\left(z_{1},z_{2}\right)}^{b}\left(\mathcal{K}_{k,r,\overline{r}}^{\left(j\right)}\right)\right\Vert _{L_{z_{1},z_{2}}^{1}\left(\mathbb{R}^{2}\times\mathbb{R}^{2}\right)}\lesssim_{a,b}\left(\frac{\lambda_{q+1}}{\overline{r}}\right)^{\left|b\right|-\left|a\right|} \label{kernel_est}
\end{align}
for $r\in\left(0,1\right)$ and $0\leq\left|a\right|,\left|b\right|\leq1$.

Because $\left|k\right|=1$, we have $k\otimes k+k^{\perp}\otimes k^{\perp}=\mathrm{Id}$.
Recall from \eqref{Rcancellation} that
\begin{align*}
\frac{R_{\ell}}{\lambda_{q+1}}+\frac{1}{2}\sum_{k\in\Omega_{i}^{q}}a_{k}^{2}\left(k^{\perp}\otimes k^{\perp}\right) & =\delta_{q+1}\mathrm{Id}\\
\iff\sum_{k\in\Omega_{i}^{q}}\frac{\lambda_{q+1}}{2}a_{k}^{2}\left(-k\otimes k\right) & =\lambda_{q+1}\left(\delta_{q+1}-\frac{1}{2}\sum_{k\in\Omega_{i}^{q}}a_{k}^{2}\right)\mathrm{Id}-R_{\ell}
\end{align*}
Then a key term in \eqref{osc_low_decomp} is
\begin{align*}
 & \sum_{i}\sum_{k\in\Omega_{i}^{q}}\frac{\lambda_{q+1}}{2}\chi_{i}^{2}\left(-k\otimes k\right)a_{k}^{2}=\sum_{i}\chi_{i}^{2}\left(\sum_{k\in\Omega_{i}^{q}}\frac{\lambda_{q+1}}{2}\left(-k\otimes k\right)a_{k}^{2}\right)\\
 & =\sum_{i} \chi_{i}^{2} \lambda_{q+1}\delta_{q+1}\mathrm{Id}-R_{\ell}-\sum_{i}\frac{\chi_{i}^{2}}{2}\sum_{k\in\Omega_{i}^{q}}a_{k}^{2}\mathrm{Id}
\end{align*}
which yields 
\begin{align*}
\divop R_{\mathrm{osc,low}} 
 & =\divop\bigg(R_{\ell}+\left(\sum_{i} \chi_{i}^{2} \lambda_{q+1}\delta_{q+1}\mathrm{Id}-R_{\ell}-\sum_{i}\frac{\chi_{i}^{2}}{2}\sum_{k\in\Omega_{i}^{q}}a_{k}^{2}\mathrm{Id}\right)\\
 &\hspace{2.2in}+\sum_{i,k}\left(\widetilde{\mathcal{Q}}_{i.k}^{\left(1\right)}+\mathcal{\widetilde{Q}}_{i.k}^{\left(2\right)}\right)\bigg)\\
 & =\divop\sum_{i,k}\left(\widetilde{\mathcal{Q}}_{i.k}^{\left(1\right)}+\mathcal{\widetilde{Q}}_{i.k}^{\left(2\right)}\right) \mod \grad C^{\infty}
\end{align*} 
where we implicitly used the fact that  $\divop\left(\sum_{i} \chi_{i}^{2}\lambda_{q+1}\delta_{q+1}\mathrm{Id}-\sum_{i}\frac{\chi_{i}^{2}}{2}\sum_{k\in\Omega_{i}^{q}}a_{k}^{2}\mathrm{Id}\right) = 0 \mod \grad C^{\infty}$. So we can pick 
\[
R_{\mathrm{osc,low}} \coloneqq \sum_{i,k}\left(\widetilde{\mathcal{Q}}_{i.k}^{\left(1\right)}+\mathcal{\widetilde{Q}}_{i.k}^{\left(2\right)}\right)
\]
and modify the pressure based on the modulo $\grad C^{\infty}$ difference.

We can now use the kernel estimates from \eqref{kernel_est},
\eqref{eq:a_psi_est}, and \eqref{eq:osc_1_ineq} to deduce 
\begin{align*}
\left\Vert R_{\mathrm{osc,low}}\right\Vert _{0} & \lesssim\sum_{i,k}\left\Vert \nabla\left(a_{k}\psi_{q+1,i,k}\right)\right\Vert _{0}\left\Vert a_{k}\psi_{q+1,i,-k}\right\Vert _{0}\\
 & \lesssim\lambda_{q}\delta_{q+1}\ll\epsilon \lambda_{q+2}\delta_{q+2},
\end{align*}

For the high-frequency part,  we argue as in Section 5.4.6 of \cite{buckmasterNonuniquenessWeakSolutions2019b}, and use \eqref{eq:a_psi_est}, \eqref{eq:w_tilde_est}, and \eqref{eq:osc_1_ineq}   to obtain
\begin{align*}
\left\Vert R_{\mathrm{osc,high}}\right\Vert _{0} & \lesssim\sum_{i,i',k,k'} \left\Vert \chi_{i}\chi_{i'} \nabla\left(a_{k}\psi_{q+1,i,k}a_{k'}\psi_{q+1,i',k'}\right)\right\Vert _{0}\\
 & \quad\quad+\left\Vert \widetilde{w}_{q+1,i,k}\right\Vert _{0}\left\Vert \nabla\left(a_{k'}\psi_{q+1,i',k'}\right)\right\Vert _{0}\\
 & \quad\quad+  \frac{1}{\lambda_{q+1}} \left\Vert \widetilde{w}_{q+1,i,k}\right\Vert_0  \left\Vert \nabla^{2}\left(a_{k'}\psi_{q+1,i',k'}\right)\right\Vert _{0}\\
 & \lesssim\lambda_{q}\delta_{q+1}+\frac{1}{\lambda_{q+1}}\lambda_{q}^{2}\delta_{q+1}\\
 &\ll \epsilon \lambda_{q+2}\delta_{q+2}.
\end{align*}
Collecting these estimates, and by arguing as above for the transport and Nash errors, we obtain \eqref{eq:R_q+1_C1}
 and \eqref{eq:R_q+1_time} for the oscillation error.  This completes the proof of Proposition~\ref{prop:iterative_proposition}.

\appendix

\section{Proof of geometric lemma}\label{lemmaproof}

We now present the proof of Lemma~\ref{lem:geometric_lemma}, which follows the same strategy as the corresponding lemma in \cite{buckmasterNonuniquenessWeakSolutions2019b}. The numerology can be explained as follows: in order to find four suitable sets of rational points on the sphere, it is convenient to use the denominator $85$, whose square can be written four ways as the sum of two squares.

Consider the rational points $k_1=(1,0)$, $k_2=(\frac{36}{85},\frac{77}{85})$, and $k_3=(\frac{36}{85},-\frac{77}{85})$ on $S^1$. One computes that we have the rank-1 decomposition
    \eqn{
    \Id=\sum_{j=1}^3a_jk_j\otimes k_j
    }
    where $a_1=4633/5929$ and $a_2=a_3=7225/11858$. The important point is not the precise values but that $a_1,a_2,a_3>0$. We define $\Omega_0^0\coloneqq\{k_1,k_2,k_3,-k_1,-k_2,-k_3\}$ and $\Omega_1^0\coloneqq\{k^\perp:k\in\Omega_0^0\}$. We perform the same operation with the rational unit vectors $m_1=(\frac{13}{85},\frac{84}{85})$, $m_2=(\frac45,\frac35)$, and $m_1=(\frac45,-\frac35)$. The corresponding rank-1 tensors generate the identity with coefficients $a_1=2023/4455$, $a_2=625/891$, and $a_3=38/45$. From these, one defines $\Omega_0^1$ and $\Omega_1^1$. As in \cite{buckmasterNonuniquenessWeakSolutions2019b}, one can perturb the positive coefficients to obtain the functions $\gamma_k$ obeying (ii) and (iii) using the implicit function theorem. One can verify directly that (iv) holds; in fact the exact lower bound is $(242/425)^{\frac12}\approx0.75$.

\section{Parameter calculus\label{sec:Parameter-calculus}}

We now list several consequences of our choice of parameters, which are used throughout the paper.  The first observation is
that
\begin{align}
\lambda_{q-1}\lambda_{q}\delta_{q-1}^{\frac{1}{2}} & \ll \tau_{c,q+1}^{-1}=\lambda_{q+1}\lambda_{q+2}\delta_{q+1}^{-\frac{1}{2}}\delta_{q+2}.\label{eq:1st_inequality}
\end{align}
Indeed, this corresponds to the condition
\begin{align*}
1-b-\beta & <  b^{2}+b^{3}-b^{2}\beta-2b^{3}\beta,
\end{align*}
which holds if and only if
\begin{align*}
\beta & < \frac{b^{3}+b^2+b-1}{2b^{3}+b^{2}-1}=1-(b-1)+O(b-1)^2.
\end{align*}
Thus \eqref{eq:1st_inequality} holds for $b$ sufficiently close to $1$ depending on $\beta$. Second, we require
\begin{equation}
{\lambda_{q-1}^{-1}}\tau_{m,q}^{-1}\delta_{q-1}^{-\frac{1}{2}}\delta_{q+1}\ll \lambda_{q+2}\delta_{q+2}\label{eq:2nd_inequality}
\end{equation}
which reduces to
\eqn{
\beta < \frac{b^2+1}{2b^2}=1-(b-1)+O(b-1)^2.
}
Next, the bound
\begin{equation}
\lambda_{q+1}^{-1}\delta_{q+1}^{\frac{1}{2}}\tau_{m,q}^{-1}\ll \lambda_{q+2}\delta_{q+2}\label{eq:3rd_ineq}
\end{equation}
is implied by \eqref{eq:2nd_inequality} upon observing that $\delta_{q-1}^\frac12\lambda_{q-1} \ll \delta_{q+1}^\frac12\lambda_{q+1}$. 

\noindent Let us also record the bound
\begin{equation}
\lambda_{q+1}^{-1}\delta_{q+1}^{\frac{1}{2}}\left(\tau_{m,q}^{-1}+\tau_{c,q+1}^{-1}\right)\ll \lambda_{q+2}\delta_{q+2}\label{eq:implication_un}
\end{equation}
which follows immediately from \eqref{eq:1st_inequality} and \eqref{eq:3rd_ineq}. The choice of $\tau_{c,q+1}$ leads to
\begin{equation}
\lambda_{q-1}^2\lambda_q^{-1}\lambda_{q+1}\tau_{c,q+1}\delta_{q-1}^{\frac{1}{2}}=\lambda_{q-1}^2\lambda_q^{-1}\lambda_{q+2}^{-1}\delta_{q+1}^{\frac{1}{2}}\delta_{q+2}^{-1}\delta_{q-1}^{\frac{1}{2}}\ll1\label{eq:1st_derived}
\end{equation}
as long as $b-1$ is so small that
\begin{align*}
\beta<\frac{b^2+b+2}{2b^2+b+1}=1-\frac{b-1}2+O(b-1)^2.
\end{align*}
Indeed, if $\beta$ and $b$ are fixed such that $\beta=1-(b-1)/2-c$, then
\eqn{
\lambda_{q-1}^2\lambda_q^{-1}\lambda_{q+1}\tau_{c,q+1}\delta_{q-1}^{\frac{1}{2}}=\lambda_q^{-4(b-1)c}\leq a^{-4(b-1)c}
}
which we can make as small as desired by choosing $a$ sufficiently large, depending only on $\beta$ and $b$. The argument is similar for the other inequalities in which we claim a ``$\ll$'' bound.

Next, the quadratic error creates the constraint
\begin{equation}
\lambda_{q}\delta_{q+1}\ll\lambda_{q+2}\delta_{q+2}\label{eq:osc_1_ineq}
\end{equation}
which follows from 
\begin{align*}
1-2b\beta & <b^{2}-2b^{2}\beta\iff\beta <\frac{b+1}{2b}=1-\frac{b-1}2+O(b-1)^2.
\end{align*}

Lastly we insist on the relation between time (frequency) scales
\begin{equation}
\tau_{c,q+1}^{-1}+\tau_{m,q}^{-1}+\lambda_{q}^{2}\lambda_{q+1}\delta_{q+1}\left(\lambda_{q-1}\delta_{q-1}^{1/2}\right)^{-1} \ll \tau_{m,q+1}^{-1}.\label{eq:M_q+1}
\end{equation}
The bound on the second term on the left is immediate from the definition. 
The first rearranges to
\eqn{
\beta\geq\frac{b^2+b+1}{2b^2+b+1}=\frac34-\frac3{16}(b-1)+O(b-1)^2.
}
While the third becomes
\[
\beta>\frac{1}{b+1}=\frac{1}{2} -\frac{b-1}{4}+O(b-1)^2
\]

With $b-1$ sufficiently small depending on the desired $\beta$, we thus conclude that the necessary constraints on the convex integration parameters are satisfied.


\end{document}